\documentclass[]{amsart}

\usepackage{cite}
\usepackage{amsmath}
\usepackage{amsfonts}
\usepackage{amssymb}
\usepackage{url}

\usepackage{algorithm}
\usepackage{algpseudocode}

\begin{document}

\title{Length-Based Attacks in Polycyclic Groups}
\author{David Garber, Delaram Kahrobaei, Ha T. Lam}

\address{David Garber, Department of Applied Mathematics, Faculty of
  Sciences, Holon Institute of Technology, 52 Golomb st., PO
  Box 305, 58102 Holon, Israel}
\email{garber@hit.ac.il}

\address{Delaram Kahrobaei, CUNY Graduate Center, PhD Program in Computer Science and NYCCT, Mathematics Department, City University of New York}
\email{dkahrobaei@gc.cuny.edu}

\address{Ha T. Lam, Department of Mathematics, CUNY Graduate Center, City University of New York}
\email{hlam@gc.cuny.edu}

\thanks{Delaram Kahrobaei is partially supported by the Office of Naval Research grant N000141210758 and also supported by PSC-CUNY grant from the CUNY research foundation, as well as the City Tech foundation.}


\begin{abstract}
The Anshel-Anshel-Goldfeld (AAG) key-exchange protocol was implemented and studied with the braid groups as its underlying platform. The length-based attack, introduced by Hughes and Tannenbaum, has been used to cryptanalyze the AAG protocol in this setting. Eick and Kahrobaei suggest to use the polycyclic groups as a possible platform for the AAG protocol.

In this paper, we apply several known variants of the length-based attack against the AAG protocol with the polycyclic group as the underlying platform. The experimental results show that, in these groups, the
implemented variants of the length-based attack are unsuccessful in the case of polycyclic groups having high Hirsch length.
This suggests that the length-based attack is insufficient to cryptanalyze the
AAG Protocol, when implemented over this type of polycyclic groups. This implies that polycyclic groups could be a potential platform for some cryptosystems based on conjugacy search problem, such as non-commutative Diffie-Hellman, El Gamal and Cramer-Shoup key-exchange protocols.

Moreover, we compare {\it for the first time} between the success rate of the different variants of the length-based attack. These experiments show that, in these groups, the memory length-based attack introduced by Garber, Kaplan, Teicher, Tsaban and Vishne does better than the other variants proposed thus far in this context.
\end{abstract}

\maketitle

\section{Introduction}
The Anshel-Anshel-Goldfeld (AAG) key-exchange protocol was introduced in 1999 \cite{AAG99}. Following its introduction, the AAG protocol was extensively studied using different groups as its underlying platform. Ko et al. \cite{KoLee00} used braid groups. Moreover, Myasnikov and Ushakov \cite{MyasnikovUshakov} studied the security of the AAG protocol with respect to several attacks on any platform groups satisfying some theoretic properties (exponentially generic free basis property).

Hughes and Tannenbaum \cite{HughesTann02} introduced the length-based attack (LBA) on the AAG protocol with its implementation in braid groups. They emphasized the importance of choosing the correct length function. Later, Garber et al. \cite{GarberLBA06} gave several realizations of this approach, particularly a length function for the braid group and experimental results suggesting that the attack fails for the parameters suggested in existing protocols. However, Garber et al. \cite{GarberProbabilistic05} also suggested an extension of the length-based attack which uses memory which succeeded in cryptanalyzing the AAG protocol. Similar attack was implemented against a system based on the Thompson group \cite{Thompson-LBA}. Most recently, Myasnikov and Ushakov \cite{MyasnikovUshakov07} analyzed the reasons behind the failure of the previous implementations of the LBA, such as the occurrence of {\it commutator-type peaks}, and gave an experimental evidence that the LBA can be modified to cryptanalyze the AAG protocol with high success rate. However, this work is again done the braid groups as the underlying platform.

Eick and Kahrobaei \cite{EickKahrobaei04} have suggested a different platform for the AAG protocol - the \textit{polycyclic group}. In polycyclic groups, the word problem can be solved efficiently \cite{Gebhardt02}, but known solutions to the conjugacy problem are much less efficient. Using experimental results, Eick and Kahrobaei showed that while the conjugacy problem can be solved within seconds using polycyclic groups with small Hirsch length, the conjugacy problem in polycyclic groups with high Hirsch length requires a much longer time for its solution.

Taking inspiration from this result, we investigate the success rate of the length-based attack on the AAG protocol, where the underlying platform is the polycyclic groups, especially those with high Hirsch length. Toward this end, we first construct polycyclic groups of high Hirsch length using a method introduced by Holt et al. \cite{Holt-handbook}. Then, we implement the different variants of the LBA presented in \cite{GarberProbabilistic05,GarberLBA06,MyasnikovUshakov07}. The experimental results that we collect suggest that the LBA is insufficient to cryptanalyze the AAG protocol, when we use the polycyclic groups with high enough Hirsch length as the underlying platform. Consequently, the polycyclic group is the {\it first} underlying platform which the LBA is insufficient for cryptanalyzing the AAG protocol on this platform, whereas the solution for the word problem is quite efficient. A suggestion for concrete parameters appears in the last section.

Moreover, we compare {\it for the first time} on any platform between the success rate of the different variants of the LBA.

As a wider application, we note that the conjugacy search problem is the basis for various cryptographic protocols besides AAG, such as the non-commutative Diffie-Hellman key-exchange \cite{KoLee00}, the non-commutative El-Gamal key-exchange \cite{KK06}, the non-abelian Cramer-Shoup key-exchange \cite{AK09} and the non-commutative digital signatures \cite{KK12}. The LBA can be applied to all these protocols; therefore, a platform group which experimental results show that the LBA is insufficient for cryptanalyzing the AAG protocol over this platform, such as polycyclic groups, can help instantiate them.

\medskip

The paper is organized as follows. In Section~\ref{AAG}, we introduce the Anshel-Anshel-Goldfeld key-exchange protocol. In Section~\ref{PolycyclicGroups}, we give a short review of polycyclic groups and the construction that we have used. In Section~\ref{LBA}, we review the length-based attack, and in Section~\ref{Results}, we present the experiments, their results and corresponding conclusions.


\section{The Anshel-Anshel-Goldfeld key-exchange protocol}
\label{AAG}
Following \cite{MyasnikovUshakov07}, we present here the Anshel-Anshel-Goldfeld key-exchange protocol (for more details, see \cite{AAG99}). As usual, we use two entities, called Alice and Bob, for presenting the two parties which plan to communicate over an insecure channel.

Let $G$ be a group with generators $g_1, \ldots, g_n$. First, Alice chooses her public set $\overline{a}=(a_1,\ldots,a_{N_1})$, where $a_i \in G$, and Bob chooses his public set $\overline{b}=(b_1,\ldots,b_{N_2})$, where $b_i \in G$. They both publish their sets. Alice then chooses her private key $A=a_{s_1}^{\varepsilon_1} \cdots a_{s_L}^{\varepsilon_L}$, where $a_{s_i} \in \overline{a}$ and $\varepsilon_i \in \{\pm1\}$. Bob also chooses his private key $B=b_{t_1}^{\delta_1} \cdots b_{t_L}^{\delta_L}$, where $b_{t_i} \in \overline{b}$ and $\delta_i \in \{\pm1\}$. Alice computes $b'_i=A^{-1}b_iA$ for all $b_i \in \overline{b}$ and sends it to Bob. Bob also computes $a'_i=B^{-1}a_iB$ for all $a_i \in \overline{a}$ and sends it to Alice. Now, the shared secret key is $K=A^{-1}B^{-1}AB$. Alice can computes this key by:
\begin{eqnarray}
K_A & = & A^{-1} (a'^{\varepsilon_1}_{s_1} \cdots a'^{\varepsilon_L}_{s_L}) = A^{-1} (B^{-1} a_{s_1} B)^{\varepsilon_1} \cdots (B^{-1} a_{s_L} B)^{\varepsilon_L} = \nonumber \\
& = & A^{-1}B^{-1} (a_{s_1}^{\varepsilon_1} \cdots a_{s_L}^{\varepsilon_L}) B = A^{-1}B^{-1}AB = K. \nonumber
\end{eqnarray}
Similarly, Bob can compute $K_B=B^{-1} (b'^{\delta_1}_{t_1} \cdots b'^{\delta_L}_{t_L}) = B^{-1}A^{-1}BA$, and then he knows the shared secret key by $K=K_B^{-1}$.

In order to find $K$, it is enough for the eavesdropper either to find $A' \in \langle a_1, \ldots, a_{N_1} \rangle$ such that $\overline{b'}=A'^{-1}\overline{b}A'$ or to find $B' \in \langle b_1, \ldots, b_{N_2} \rangle$ such that $\overline{a'}=B'^{-1}\overline{a}B'$ (an incompatible sufficient condition can be found in \cite{KTV}). Thus, the security of the AAG protocol is based on the assumption that the subgroup-restricted simultaneous conjugacy search problem is hard.


\section{Polycyclic groups}
\label{PolycyclicGroups}
In this section, we give a short review for polycyclic groups and describe the construction of polycyclic groups of high Hirsch length. For more details, see \cite{Holt-handbook}.

\subsection{The polycyclic presentation} Recall that $G$ is a \textit{polycyclic group} if it has a polycyclic series, i.e., a subnormal series
$$G=G_1 \rhd G_2 \rhd \cdots \rhd G_{n+1}=\{1\},$$
with non-trivial cyclic factors. The {\it polycyclic generating sequence} of $G$ is the $n$-tuple $(g_1,\ldots, g_n)$, such that $G_i=\langle g_i,G_{i+1} \rangle$ for $1\leq i \leq n$.

Any polycyclic group has a finite presentation of the form:
$$\langle g_1,\ldots,g_n \mid g_j^{g_i}=w_{ij}, \; g_j^{g_i^{-1}}=v_{ij}, \; g_k^{r_k}=u_{kk} \; \text{for} \; 1\leq i<j\leq n \; \text{and} \; k \in I \rangle $$
where $w_{ij},v_{ij},u_{kk}$ are words in the generators $g_{i+1},\ldots ,g_n$ and $I$ is the set of indices $i \in \{1,\ldots,n\}$ such that $r_i=[G_i:G_{i+1}]$ is finite. Here $a^b$ stands for $b^{-1}ab$.

It is known by induction that each element of $G$ defined by this presentation can be uniquely written as $g=g_1^{e_1} \cdots g_n^{e_n}$ where $e_i \in \mathbb{Z}$ for $1\leq i \leq n$, and $0 \leq e_i < r_i$ for $i \in I$. We call $g=g_1^{e_1} \cdots g_n^{e_n}$ the \textit{normal form} of an element in $G$. If every element in the group can be uniquely presented in the normal form, then the polycyclic presentation is called \textit{consistent}. Note that every polycyclic group has a consistent polycyclic presentation \cite{Holt-handbook}.

The \textit{Hirsch length} of a polycyclic group is the number of indices $i$ such that $r_i=[G_i:G_{i+1}]$ is infinite. This number is invariant of the chosen polycyclic sequence.

\subsection{Constructing polycyclic groups using number fields}\label{3.2}
There are several ways for constructing polycyclic groups. For the purpose of this paper, we construct polycyclic groups by semidirect products of the maximal order and the unit group of a number field. This construction follows \cite{Holt-handbook}.

Let $f(x) \in \mathbb{Z}[x]$ be an irreducible polynomial. The polynomial $f$ defines a field extension $F$ over $\mathbb{Q}$. The \textit{maximal order} or the \textit{ring of integers} $O_F$ of the number field $F$ is the set of algebraic integers in $F$:
$$O_F=\{a \in F \mid \text{ there exists a monic polynomial } f_a(x) \in \mathbb{Z}[x] \text{ such that } f_a(a)=0 \}.$$
The \textit{unit group} of $F$ is:
$$U_F=\{ a \in O_F \mid a \neq 0 \text{ and } a^{-1} \in O_F \}.$$

For constructing the polycyclic group by the maximal order and the unit group of a number field $F$ where $[F:\mathbb{Q}]=n$, we recall two results. First, the maximal order $O_F$ forms a ring whose additive group is isomorphic to $\mathbb{Z}^n$ \cite{StewartTall}. Second, Dirichlet's unit theorem states that given $n=s+2t$, where $s$ and $2t$ are the numbers of real field monomorphisms $F \rightarrow \mathbb{R}$ and complex field monomorphisms $F \rightarrow \mathbb{C}$ respectively, then the unit group $U_F$ is a finitely-generated abelian group of the form $U_F \cong \mathbb{Z}^{s+t-1} \times \mathbb{Z}_m$ for some even $m$ \cite{StewartTall}. Here, we use the fact that the unit group is a finitely-generated abelian group and hence $U_F$ is also polycyclic.

Let $G$ be a group and $N \trianglelefteq G$, it is easy to see that if $N$ and $G/N$ are both polycyclic, then the group $G$ is also polycyclic by putting together the polycyclic series of $N$ and the series induced by the polycyclic series of $G/N$. Since the above results guaranteed that the maximal order is a polycyclic group and the unit group, which is isomorphic to $G/O_F$, is also polycyclic, the group $G=O_F \rtimes U_F$ is polycyclic. The action which defines the semidirect product is a multiplication from the right of $U_F$ on $O_F$.

If $N \trianglelefteq G$ ,the Hirsch length of a polycyclic group $G$ is $h(G)=h(N)+h(G/N)$; in our case, $h(G)=h(O_F) + h(U_F)$, where $h(O_F)$ is $n$, which is the degree of the generating polynomial $f$. Hence, for constructing a polycyclic group of high Hirsch length, we have to find an irreducible polynomial of high enough degree, and then the polycyclic group constructed by the above method will have Hirsch length larger than the degree of the polynomial.

\subsection{Polycyclic groups as platform groups for the AAG protocol}
Polycyclic groups are suitable as platform groups for the AAG protocol for several reasons. First, the word problem can be solved efficiently using the collection algorithm \cite{Gebhardt02}, see also \cite{EickKahrobaei04}. Second, the conjugacy search problem has no efficient solution in general polycyclic groups. This assessment is due to Eick and Kahrobaei \cite{EickKahrobaei04}, using the following experiment: let $K=\mathbb{Q}[x]/(f_w)$ be an algebraic number field for a cyclotomic polynomial $f_w$, where $w$ is a primitive $r$-th root of unity. Let $G(w)=O \rtimes U$, where $O$ is the maximal order and $U$ the unit group of $K$, $r$ the order of $w$ and $h(G(w))$ the Hirsch length. The average time used for 100 applications of the collection algorithm on random words and the average time used for 100 applications of the conjugacy algorithm on random conjugates are:

\medskip

\begin{center}
\begin{tabular}{|c|c|c|c|}
\hline
 r & h(G(w)) & Collection & Conjugation \\
\hline 3 & 2 & 0.00 seconds & 9.96 seconds \\
 4 & 2 & 0.00 seconds & 9.37 seconds \\
 7 & 6 & 0.01 seconds & 10.16 seconds \\
 11 & 14 & 0.05 seconds & $>$ 100 hours \\
\hline
\end{tabular}
\end{center}

\medskip

We can see that the collection algorithm works very fast even for polycyclic groups of high Hirsch length, and therefore the word problem has an efficient solution. On the other hand, the solution to the conjugacy problem is not efficient for polycyclic groups having high Hirsch length.


\section{The length-based attack}
\label{LBA}

The length-based attack (LBA) is a probabilistic attack against the conjugacy search problem in general, and against the AAG protocol in particular, with the goal of finding Alice's (or Bob's) private key. It is based on the idea that a conjugation of the correct element should decrease the length of the captured package. Using the notations of Section ~\ref{AAG}, the captured package is $\overline{b'}=(b'_1, \ldots , b'_{N_2})$, where $b'_i=A^{-1}b_iA$. If we conjugate $\overline{b'}$ with elements from the group $\langle a_1, \ldots, a_{N_1} \rangle$ and the length of the resulting tuple has been decreased, then we have found a candidate for the conjugating factor. The process of conjugation is then repeated with the decreased-length tuple until a longer candidate for the conjugating factor is found. The process ends when the conjugated captured package is the same as $\overline{b}=(b_1, \ldots , b_{N_2})$, which is known. Then, the conjugate can be recovered by reversing the sequence of conjugating factors. For more details on the LBA, see \cite{GarberProbabilistic05,GarberLBA06,HT10,Myasnikov-book, MyasnikovUshakov07}.

\subsection{Variants of the LBA}
\label{LBAVariations}
In \cite{GarberProbabilistic05,GarberLBA06,MyasnikovUshakov07,Thompson-LBA}, several variants of the LBA are presented. Here, we give four variants of the LBA that we implemented against the AAG protocol having the polycyclic group as its underlying platform. In all these variants, the following input and output are expected:
\begin{itemize}
\item \textsc{Input}: $\overline{a}=(a_1, \ldots , a_{N_1})$, $\overline{b}=(b_1, \ldots , b_{N_2})$ and $\overline{b'}=(b'_1, \ldots , b'_{N_2})$, such that $b'_i=b_i^A$ for $i=1,\ldots,N_2$.
\item \textsc{Output}: An element $A' \in \langle a_1, \ldots, a_{N_1} \rangle$ such that $b'_i=b_i^{A'}$ for $i=1,\ldots,N_2$, or FAIL if the algorithm cannot find such $A'$.
\end{itemize}

We will use the following notation: if $\overline{c}=(c_1,\ldots,c_k)$, then its \textit{total length} $|\overline{c}|$ is $\sum_{i=1}^k |c_i|$ (the length of $c_i$, $|c_i|$, will be discussed in Section \ref{sec_length}).

\subsubsection{LBA with backtracking}
The most straight-forward variant of LBA (Algorithm ~\ref{algo1}) conjugates $\overline{b'}$ directly with $a_i^{\pm 1} \in \{a_1, \ldots, a_{N_1}\}$. This is termed ``LBA with backtracking'' by Myasnikov and Ushakov \cite{MyasnikovUshakov07}.

\begin{algorithm}
\caption{LBA with backtracking}
\label{algo1}
\begin{algorithmic}[1]
	\State Initialize $S=\{(\overline{b'},\rm{id_G})\}$.
	\While{$S \neq \emptyset$}
		\State Choose $(\overline{c},x) \in S$ such that $|\overline{c}|$ is minimal. Remove $(\overline{c},x)$.
		\For{$i=1,\ldots, N_1$ and $\varepsilon=\pm 1$}
			\State Compute $\delta_{i,\varepsilon}=|\overline{c}|-\left|\overline{c}^{a_i^\varepsilon}\right|$.
			\State \textbf{if} $\overline{c}^{a_i^\varepsilon}=\overline{b}$ \textbf{then} output inverse of $xa_i^\varepsilon$ and stop.
			\If{$\delta_{i,\varepsilon}>0$} \Comment{length has been decreased}
				\State Add $\left(\overline{c}^{a_i^\varepsilon},xa_i^\varepsilon\right)$ to $S$.
			\EndIf
		\EndFor
	\EndWhile
	\State Otherwise, output FAIL. \Comment no more elements to conjugate
\end{algorithmic}
\end{algorithm}

\subsubsection{LBA with a dynamic set}
Through analysis, Myasnikov and Ushakov \cite{MyasnikovUshakov07} concluded that different types of peaks make LBA unsuccessful. To overcome this, they suggested a new version of the algorithm, which they termed ``LBA with a dynamic set''. Here, if a generator $a_i$ causes a length reduction, only the conjugates and products involving $a_i$ are added to the dynamic set. On the other hand, if no generator causes a length reduction, all conjugates and two generators products are added. Their experimental results suggest that this algorithm works especially well in the case of keys composed from long generators, but it is not worse than the naive  algorithm in the other cases. The algorithm presented here is a modified version of their algorithm, which we implemented to attack the AAG protocol having the polycyclic group as its underlying platform.

\begin{algorithm}
\caption{LBA with a dynamic set}
\label{algo5}
\begin{algorithmic}[1]
	\State Initialize $S=\{(\overline{b'},\rm{id_G})\}$.
	\While{$S \neq \emptyset$}
		\State Choose $(\overline{c},x) \in S$ such that $|\overline{c}|$ is minimal. Remove $(\overline{c},x)$.
		\For{$i=1,\ldots, N_1$ and $\varepsilon=\pm 1$}
			\State Compute $\delta_{i,\varepsilon}=|\overline{c}|-\left|\overline{c}^{a_i^\varepsilon}\right|$
		\EndFor
		\If{$\delta_{i,\varepsilon} \leq 0$ for all $i$}
			\State Define $\overline{a}_{\rm{ext}} = \overline{a} \cup \{x_i x_j x_i^{-1}, x_i x_j, x_i^2 \mid x_i, x_j \in \overline{a}^{\pm 1}, i \neq j \}$.
		\Else $\,$ Define $\overline{a}_{\rm{ext}} = \overline{a} \cup \{x_j x_m x_j^{-1}, x_m x_j, x_j x_m, x_m^2 \mid x_j \in \overline{a}^{\pm 1}, m \neq j \}$ where $x_m$ such that $\delta_m = \text{max}\{\delta_{i,\varepsilon} \mid i=1,\ldots, N_1 \}$.
		\EndIf
		\ForAll{ $w \in \overline{a}_{\rm{ext}}$ }
			\State Compute $\delta_w=|\overline{c}|-|\overline{c}^w|$.
		\EndFor
		\State \textbf{if} $\overline{c}^{w}=\overline{b}$ \textbf{then} output inverse of $xw$ and stop.
		\If{$\delta_w>0$} \Comment{length has been decreased}
			\State Add $(\overline{c}^w,xw)$ to $S$.
		\EndIf
	\EndWhile
	\State Otherwise, output FAIL. \Comment no more elements to conjugate
\end{algorithmic}
\end{algorithm}

\subsubsection{Memory-LBA}
Another variant, presented in \cite{GarberProbabilistic05}, is also considered. In this variant, we allocate an array $S$ of a fixed size $M$. The array $S$ holds $M$ tuples every round. In every round, all elements of $S$ are conjugated, but only the $M$ smallest conjugated tuples (with respect to their length) are inserted back into $S$. For the halting condition, we use a predefined time-out.

\begin{algorithm}
\caption{Memory-LBA}
\label{algoMem2}
\begin{algorithmic}[1]
	\State Initialize $S=\{(|\overline{b'}|,\overline{b'},\rm{id_G})\}$.
	\While{not time-out}
		\For{$(|\overline{c}|,\overline{c},x) \in S$}
			\State Remove $(|\overline{c}|,\overline{c},x)$ from $S$.
			\State Compute $\overline{c}^{a_i^\varepsilon}$ for all $i \in \{1 \ldots N_1\}$ and $\varepsilon \in \{ \pm 1 \}$.
			\State \textbf{if} $\overline{c}^{a_i^\varepsilon}=\overline{b}$ \textbf{then} output inverse of $xa_i^\varepsilon$ and stop.
			\State Save $\left(\left|\overline{c}^{a_i^\varepsilon}\right|,\overline{c}^{a_i^\varepsilon},xa_i^\varepsilon\right)$ in $S'$.
		\EndFor
		\State After finished all conjugations, sort $S'$ by the first element of every tuple
		\State Copy the smallest $M$ elements into $S$ and delete the rest of $S'$
	\EndWhile
	\State Otherwise, output FAIL.
\end{algorithmic}
\end{algorithm}

\subsubsection{LBA* (with memory)}
We present a different variant of memory-LBA which is again based on a fixed-size array allocated for the algorithm. 
Here, $S$ holds $M$ tuples every round and is sorted by the first element (with respect to the length of conjugated element) of each tuple. In every round, only the smallest element of $S$ is removed and conjugated by all the generators and their inverses. The conjugated tuples are inserted back into $S$ depending on whether there is a free place in $S$. If there is no more places in $S$, and the conjugated tuple is smaller than the largest element in $S$, swap them and re-sort $S$. Since $S$ is always kept sorted, any operation to find the ``smallest element'' costs constant time. As in the previous variant, we use a predefined time-out as the halting condition.

The name LBA* comes from the general idea of {\it A* search algorithm} \cite{HNR}, which uses a best-first search (as we are doing here - taking the smallest element of $S$ and conjugated it). We should note that a very similar algorithm was independently introduced by Tsaban \cite{TsabanSlides}, and the difference between the two variants is that our variant starts the search from $\overline{b'}$, while Tsaban's variant starts the search from both directions: $\overline{b'}$ and $\overline{b'}$ (using the idea of ``meet in the middle'').

\begin{algorithm}
\caption{LBA* (with memory)}
\label{algoMem1}
\begin{algorithmic}[1]
	\State Initialize $S=\{(|\overline{b'}|,\overline{b'},\rm{id_G})\}$.
	\While{not time-out}
		\State Choose $(|\overline{c}|,\overline{c},x) \in S$ such that $|\overline{c}|$ is minimal. Remove $(|\overline{c}|,\overline{c},x)$.
		\For{$i=1,\ldots, N_1$ and $\varepsilon=\pm 1$}
			\State Compute $\overline{c}^{a_i^\varepsilon}$.
			\State \textbf{if} $\overline{c}^{a_i^\varepsilon}=\overline{b}$ \textbf{then} output inverse of $xa_i^\varepsilon$ and stop.
			\If{$\rm{Size}(S) < M$}
				\State Add $\left(\left|\overline{c}^{a_i^\varepsilon}\right|,\overline{c}^{a_i^\varepsilon},xa_i^\varepsilon\right)$ to $S$ and sort $S$ by first element of every tuple.
			\Else \Comment no more space in S
				\State \textbf{if} $\left|\overline{c}^{a_i^\varepsilon}\right|$ is smaller than first element of all tuples in S \textbf{then} swap them
			\EndIf
		\EndFor
	\EndWhile
	\State Otherwise, output FAIL. \Comment no more elements to conjugate
\end{algorithmic}
\end{algorithm}

\subsection{The length function}\label{sec_length}
In the implementation of the LBA, the choice of the length function is important (see \cite{GarberProbabilistic05,HT10}). In our case, the length of a word is chosen to be the sum of the absolute values of the exponents in its normal form. We choose this function because the experimental results presented below show that it satisfies the requirement $\ell(a^{-1}ba) \gg \ell(b)$ (as needed for a length function used for LBA).

The first step of the experiments is the construction of  a polycyclic group $G$ of a given Hirsch length $h(G)$, following the construction in Sections \ref{3.2} and \ref{PolyConstruction}. Then, an element $b$ of length between 10 and 13 is randomly chosen; we choose elements of this length for consistency with the LBA parameters. Another random element $a$ satisfying the same length interval is chosen and $b^a$ is computed, and finally, we compute $|b^a|-|b|$. We performed 100 tests for each group and the average difference is recorded.

\medskip

\begin{center}
\begin{tabular}{|c|c|c|}
\hline
 Polynomial & h(G) & Average difference \\
\hline
 $x^2-x-1$ & 3 & 79.92 \\
 $x^5-x^3-1$ & 7 & 80.17\\
 $x^{11}-x^3-1$ & 16 & 44.93 \\
\hline
\end{tabular}
\end{center}

\medskip

As we can see, the average difference is large; specifically $|b^a|-|b|$ is significantly larger than $|a|$, indicating that the condition $\ell(a^{-1}ba) \gg \ell(b)$ is indeed satisfied.


\section{Experimental results}
\label{Results}

Our goal is to apply the LBA on the AAG protocol having the polycyclic group as its underlying platform. To that end, we implemented the four variants of the LBA presented in Section ~\ref{LBA} and performed experiments on several polycyclic groups having different Hirsch lengths.

\subsection{Implementation details}
\label{PolyConstruction} Each polycyclic group is constructed by choosing an irreducible polynomial $f$ over $\mathbb{Z}$, thus $f$ defines an algebraic field $F$ over $\mathbb{Q}$. Let $O_F$ be its maximal order and $U_F$ be its unit group, thus $O_F \rtimes U_F$ is the desired polycyclic group. This construction follows \cite{Holt-handbook} and is a part of the Polycyclic package of GAP \cite{GAPPolycyclic}.

A random element $a_i$, for Alice's public set, or $b_i$, for Bob's public set, is generated by taking either some random generators of the group or their inverses and multiplying them together, while maintaining that the length of the element is between a predefined minimum and maximum. By this method, we take control over the length of the element.

Alice's private key $A$ is generated by taking a fixed number of random elements in $\overline{a}=(a_1,\ldots,a_{N_1})$ and multiplying them together. Here we forgo control over length to preserve interesting cases of conjugations actually decreasing the length of $b_i$, such as a commutator-type peak. The way for choosing the keys is similar to what has been used in \cite{MyasnikovUshakov07}. This way also reflects the characterization of the polycyclic group.

\subsection{Results}
We performed several sets of tests, all of which were run on an Intel Core I7 quad-core 2.0GHz computer with 12GB of RAM, running Ubuntu Version 12.04 with GAP Version 4.5 and 10GB of memory allocation. In all these tests, the polycyclic group $G$ having Hirsch length $h(G)$ is constructed by the above method using polynomial $f$. The size of Alice's and Bob's public sets are both $N_1=N_2=20$.

\subsubsection{The effect of the Hirsch length}
In the first set of tests, the length of each random element $a_i$ or $b_i$ is in the interval $[L_1,L_2]=[10,13]$ and Alice's private key is the product of $L=5$ random elements in Alice's public set. The time for each batch of 100 tests are recorded together with its success rate. In each case, a time-out of 60 minutes is enforced for each test. The following results are obtained by LBA with a dynamic set:

\medskip

\begin{center}
\begin{tabular}{|c|c|c|c|}
\hline
 Polynomial & h(G) & Time & Success rate of \\
            &      &      & LBA with a dynamic set \\
\hline
 $x^2-x-1$ & 3 ~ & 0.20 hours & 100\% \\
 $x^5-x^3-1$ & 7 & 76.87 hours & 35\% \\
 $x^7-x^3-1$ & 10 & 94.43 hours & 8\% \\
 $x^9-7x^3-1$ & 14 & 95.18 hours & 5\% \\
 $x^{11}-x^3-1$ & 16 & 95.05 hours & 5\% \\
\hline
\end{tabular}
\end{center}

\medskip

From this table, we can see that with a small Hirsch length, the LBA cryptanalyzes the AAG protocol easily with high success rate. However, as the Hirsch length is increased to 7, the success rate decreases. In polycyclic groups with higher Hirsch lengths, we can see the effect of the time-out more prominently as the total time did not increase much more, but the success rate is dropped to 5\%. Although a success rate of 5\% is not negligible, note that we use a very small value for $L$. Based on the current experimental results, we expect that increasing the value of $L$ will reduce the success rate to 0\%.

\subsubsection{The effect of the key length}
In the second set of tests, we vary the number of elements $L$ that compose Alice's private key. Myasnikov and Ushakov \cite{MyasnikovUshakov07} suggested that the LBA with a dynamic set has a high success rate with long generators, i.e. random elements have longer length $[L_1,L_2]$. Therefore, we also vary the length of random elements according to the parameters in \cite{MyasnikovUshakov07}.

The following results are obtained by LBA with a dynamic set, with a time-out of 30 minutes:

\medskip

\begin{center}
\begin{tabular}{|c|c|c|c|c|c|c|}
\hline
	Polynomial & h(G) & [10,13] & \multicolumn{2}{|c|}{[20,23]} & [40,43] \\
\hline
					&	&$L=10$	& $L=10$&$L=20$	& $L=50$ \\
\hline
	$x^7-x^3-1$ 		& 10& 2\%	& 0\%	& 0\%	& 0\% \\
	$x^9-7x^3-1$		& 14& 0\%	& 0\%	& 0\%	& 0\% \\
	$x^{11}-3x^3-1$ 	& 17& 0\%	& 0\%	& 0\%	& 0\% \\
\hline
\end{tabular}
\end{center}

\medskip

The results of this set of tests indicate that just by increasing the number of generators of Alice's private key from 5 (as in the previous set of tests) to 10, the LBA already fails with polycyclic groups having Hirsch length as small as 10.

\subsubsection{Comparing the four variants of the LBA}
In this paper, we compare the success rate of the four variants of the LBA {\it for the first time} on any platform. For comparing the success rate of the four variants of the LBA, we purposely choose the value of the test parameters to be very small in this set of tests. They are as follows: $N_1=N_2=20$, $[L_1,L_2]=[5,8]$, $L=5$, there is a time-out of 30 minutes and a memory of size $M=500$. The polynomial used is $f=x^3-x-1$, constructing a polycyclic group of Hirsch length 4.

\medskip

\begin{center}
\begin{tabular}{|l|c|c|c|}
\hline
Algorithm & Time & Success rate \\
\hline
LBA with backtracking & 0.57 hours & 58\% \\
LBA with a dynamic set & 37.35 hours & 95\% \\
Memory-LBA (with memory $M=500$) & 4.01 hours & 92\% \\
LBA* (with memory $M=500$) & 32.00 hours & 36\% \\
\hline
\end{tabular}
\end{center}

\medskip

Algorithm LBA with a dynamic set gives the best success rate but took much longer than Algorithm Memory-LBA which gives a similar success rate in much shorter time. We conclude that with a sufficient size of memory, Algorithm Memory-LBA is the best variant of the LBA.

\subsubsection{Using the four variants of the LBA on our test parameters}
In the fourth set of tests, we want to see the effect of the four different variants of the LBA presented in Section ~\ref{LBAVariations} applied to our test parameters. Therefore, we keep the following parameters for all the algorithms: the length of each random element is in the interval $[L_1,L_2]=[10,13]$, Alice's private key is the product of 10 elements and the length of both public sets are $N_1=N_2=20$. There is a time-out of 30 minutes per test and in the case of the two memory variants of the LBA, Algorithm Memory-LBA and Algorithm LBA*, a memory of size $M=1000$ is used. The same polycyclic group $G$ having Hirsch length 14 constructed by the polynomial $x^9-7x^3-1$ is used for all the variants of the LBA.

\medskip

\begin{center}
\begin{tabular}{|l|c|c|}
\hline
 Algorithm & Time & Success rate \\
\hline
LBA with backtracking & 48.68 hours & 0\% \\
LBA with a dynamic set & 50.04 hours & 0\% \\
Memory-LBA (with memory $M=1000$) & 49.35 hours & 3\% \\
LBA* (with memory $M=1000$) & 50.00 hours & 0\% \\
\hline
\end{tabular}
\end{center}

\medskip

As we can see, Memory-LBA algorithm has the best performance in this set of parameters, but even then, it has only 3\% success rate. To further test Memory-LBA algorithm, we ran another set of tests where we increase the length of random elements to $[L_1,L_2]=[20,23]$ and increase the number of factors of the private key to $L=20$. To give it a chance of success, we increase the size of the memory $M$ to 40,000. The result is 0\% success rate.

\subsubsection{The effect of increasing the time-out}
Since it is possible that the time-out of 30 minutes for each test is too short, we ran another set of tests, where the time-out is 4 hours for each test. Memory-LBA algorithm showed the most promise, so we chose it with the following parameters: the length of random elements is in the interval $[L_1,L_2]=[20,23]$, the number of factors of the private key is $L=20$ and the size of the memory $M$ is 1000. The polynomial used is $x^9-7x^3-1$ producing a polycyclic group of Hirsch length 14. Due to the long time-out, we performed only 50 tests. We still get 0\% success rate.

\medskip

Based on the above experimental results, we conclude that the LBA is insufficient for cryptanalyzing the polycyclic groups of high enough Hirsch lengths. One can suggest the following parameters: $h(G)=16, L=20$ and $[L_1,L_2]=[20,23]$ for achieving an AAG protocol based on the polycyclic group, which the known variants of the LBA have 0\% success rate for cryptanalyzing this protocol.

\subsubsection{Additional experimental results concerning LBA with a dynamic set algorithm}
Here, we present some additional experimental results for LBA with a dynamic set. The time-out for each test is $1$ hour. The polynomials used are $f$ and $h(G)$ is the Hirsch length of the corresponding polycyclic group. The size of Alice's and Bob's public sets are $N_1,N_2$ respectively. Each random element $a_i$ or $b_i$ has length in $[L_1,L_2]$ and Alice's private key is the product of $L=5$ random elements in Alice's public set. The success rate of a batch of 100 tests is recorded.

\medskip

\begin{center}
\begin{tabular}{|c|c|c|c|c|}
\hline
	Polynomial & h(G) & \multicolumn{2}{|c|}{$N_1=N_2=5$} & $N_1=N_2=20$ \\
\hline
	& & [5,8] & [15,18] & [10,13] \\
\hline
	$x-1$ 			& 1	& 98\%	& 		& 98\%	\\
	$x^2-x-1$ 		& 3 & 98\%	& 96\%	& 100\%	\\
	$x^3-x-1$ 		& 4 & 95\%	& 		& 100\%	\\
	$x^5-x^3-1$ 		& 7 & 		& 		& 35\% 	\\
	$x^7-x^3-1$ 		& 10& 		& 		& 8\%	\\
	$x^9-7x^3-1$		& 14& 		& 		& 5\%	\\
	$x^{11}-x^3-1$ 	& 16& 59\%	& 53\%	& 5\%	\\
\hline
\end{tabular}
\end{center}

\section*{Acknowledgements}
We would like to thank an anonymous referees for many useful suggestions, which were implemented in the text. 


\bibliographystyle{plain}
\bibliography{PolycyclicLBApaper-arxiv}

\end{document}